\documentclass[11pt]{amsart}

\usepackage[T1]{fontenc}
\usepackage[ansinew]{inputenc}

\usepackage{euler}

\usepackage{amsmath}  
\usepackage{amsthm}   
\usepackage{latexsym}
\usepackage{amsfonts}
\usepackage{txfonts}
\usepackage{mathrsfs}

\usepackage{graphicx}
\usepackage{psfrag}

\numberwithin{figure}{section}

\numberwithin{equation}{section}

\begin{document}

\newtheorem{teor}{Theorem}[section]
\newtheorem{defn}[teor]{Definition}
\newtheorem{prop}[teor]{Proposition}
\newtheorem{exmp}[teor]{Example}
\newtheorem{lemm}[teor]{Lemma}
\newtheorem{seur}[teor]{Corollary}
\newtheorem{exmps}[teor]{Examples}

\title[Solutions of Aronsson equation near isolated points]
{Solutions of Aronsson equation near isolated points}

\author{Vesa Julin}

\address{Department of Mathematics and Statistics,
P.O.Box 35, FIN-40014 University of Jyv\"askyl\"a, Finland}
\email{vvjulin@cc.jyu.fi}


\date{\today}

\begin{abstract} When  \( n \geq 2 \), we show that for a non-negative solution of the Aronsson equation \( \mathscr{A}_H (u) = D_x H(Du(x)) \, \cdot \, H_p(Du(x)) = 0 \) an isolated singularity \(x_0 \) is either a removable singularity or \( u(x) = b + C_{k}^{H}(x-x_0)  + o(|x - x_0|) \) ( or \( u(x) = b - C_{k}^{\hat{H}}(x-x_0)  + o(|x - x_0|) \) ) for some \( k > 0 \) and \( b \in \mathbb{R}\) . Here \(C_{k}^{H}\) and  \( C_{k}^{\hat{H}} \) are general cone functions. This generalizes the asymptotic behavior theory for infinity harmonic functions by Savin, Wang and Yu \cite{SaWaYu}. The Hamiltonian \( H \in C^2(\mathbb{R}^n) \) is assumed to be non-negative and uniformly convex.
\end{abstract}
\maketitle

\section{Introduction}
\vspace*{4mm}
Variational problems for \(L^{\infty}\)-functionals 
\[
 F(u,\Omega) = \underset{x \in \Omega}{ \text{ess sup}} \,  H (Du(x), u(x), x ) \, , 
\]
where \( u \in W^{1, \infty}(\Omega) \), were first studied by Aronsson in 1960's (see \cite{Aronsson1} - \cite{Aronsson2} ). He formally derived the corresponding Euler-Lagrange equation
\begin{equation}
 \label{eulerarons}
 D_x H(Du(x),u(x),x) \, \cdot \, H_p(Du(x), u(x),x) = 0
\end{equation}
for minimizers of \( F( \cdot,\Omega)\) which he called \emph{absolute minimizers}. Here \(H_p \) denotes the gradient of \( H(p,s,x) \) with respect to the first variable and \(D_x H(Du(x),u(x),x) \) is the gradient of the map \(x \mapsto H(Du(x),u(x),x) \). Nowadays the equation (\ref{eulerarons}) is known as Aronsson equation. Aronsson also noticed the importance of the special case
\( H(p) = \frac{1}{2}|p|^2 \). Then (\ref{eulerarons}) has the form
\begin{equation}
\label{aaretonlaplace}
  \triangle_{\infty} u(x) = D^2 u(x) D u(x) \cdot D u(x) = 0
\end{equation}
and the solutions are called infinity harmonic functions. However Aronsson also noted that even in this special case these solutions might not be smooth.

A major step forward was taken by Jensen in \cite{jensen}, where he suggested that the equation (\ref{aaretonlaplace}) should be interpreted in the viscosity sence. He proved the equivalence of infinity harmonic functions and absolute minimizers of the functional \( F(u,\Omega) = || Du(x)||_{L^{\infty}(\Omega)}\). He also proved the existence and uniqueness of an absolute minimizer with given continuous boundary values. Another major step was taken by Crandall, Evans and Gariepy \cite{cran-evans-gar}, who introduced the concept of comparison with cones, which turned out to be a very useful tool. E.g. by using this characterization Savin and Evans (in \cite{savin} and \cite{C-1-alfa}  ) were able to prove the local \( C^{1,\alpha} \) regularity for infinity harmonic functions in two dimensions. One of the most interesting open questions in the field is whether the same regularity result is true also in higher dimensions. For those readers who are not familiar with infinity harmonic functions we refer to an excellent survey by Aronsson, Crandall and Juutinen \cite{ArGrJu} which covers more or less the whole basic theory. 

In recent years there has been an increasing interest in slightly more general Aronsson equation 
\begin{equation}
 \label{aronsson2}
 \mathscr{A}_H (u) = D_x H(Du(x)) \, \cdot \, H_p(Du(x)) = 0 ,
\end{equation}
where the Hamiltonian \( H \) is usually assumed to have some convexity properties. Gariepy, Wang and Yu \cite{generalcomparisonGWY} introduced the concept of comparison with general cones and proved the equivalence between solutions of (\ref{aronsson2}) and absolute minimizers of the functional \(F(u,\Omega) = || H(Du(x))||_{L^{\infty}(\Omega)} \). By using this general comparison principle Wang and Yu \cite{WaYu} generalized the original regularity result of Savin by showing that in dimension two solutions of (\ref{aronsson2}) are locally of class \(C^1 \), when the Hamiltonian  \( H \) is assumed to be uniformly convex. 

In this paper we study solutions of (\ref{aronsson2}) near an isolated point. The Hamiltionian is assumed to be a non-negative and uniformly convex \( C^2 \) function. Our first theorem is a generalization of a result of Bhattacharya \cite{Bhattacharya}. 
\begin{teor}
 \label{bhatta2}
Let \( u \in C (B_R(x_0) \backslash  \{x_0 \})  \) be a non-negative viscosity solution of
\[
 \mathscr{A}_H (u)= 0 \qquad \text{in} \, B_R(x_0) \backslash  \{x_0 \}  .
\]
Then the limit
\[
 \lim_{x \to x_0} u(x)
\]
exists. 
\end{teor}
One can not drop the assumption of non-negativity, since by Aronsson \cite{Aronsson3} there are infinity harmonic functions which are unbounded near an isolated point. 

Our second main result deals with the asymptotic behavior of solutions of (\ref{aronsson2}) and is strongly motivated by Savin, Wang and Yu \cite{SaWaYu}. With a help of Theorem \ref{bhatta2} we show that a solution of Aronsson equation behaves as a general cone near a singular point. 
\begin{teor}
\label{teoria1}
 Suppose \( n \geq 2 \). Let a function \( u \geq 0 \) be a viscosity solution of the Aronsson equation \( \mathscr{A}_H (u) = 0\) in \( B_r(x_0) \backslash \{x_0\}\). Then one of the two alternatives holds: 
\begin{itemize}
\item[(i)] \( u\) is a viscosity solution of \( \mathscr{A}_H (u) = 0\) in the whole ball \(B_r(x_0)\) 

\item[(ii)] either 
\[
 u(x) = b + C_{k}^{H}(x-x_0)  + o(|x - x_0|)
\]
or
\[
  u(x) = b - C_{k}^{\hat{H}}(x-x_0)  + o(|x - x_0|) ,
\]
for some \( k > 0 \) and \(b \geq 0 \).
\end{itemize}
\end{teor}
Above \( C_{k}^{H}\), \( C_{k}^{\hat{H}} \) are general cone functions, which will be defined later in Section 2. This is a far deeper result than Theorem \ref{bhatta2} and the proof is more challenging too. Compared to Theorem 1.1 in \cite{SaWaYu} where \( u \) is infinity harmonic, the proof needs new technical methods and more careful geometric arguments.

With a help of Theorem \ref{teoria1} and some estimates from \cite{WaYu} we derive a Corollary which gives us a family of nonclassical solutions of (\ref{aronsson2}).

\begin{seur}
\label{corollaarimain}
 Suppose \( \Omega \subset \mathbb{R}^2 \) is a bounded domain and \( x_0 \in \Omega \). Let a function \( u \in C(\overline{\Omega} )\) be a solution of the Aronsson equation  \( \mathscr{A}_H (u) = 0\) in \( \Omega \backslash \{x_0 \} \) which takes boundary values \( u = 1 \) in \( \partial \Omega \) and \( u(x_0) = 0 \). Then \( u \in C^2( \Omega \backslash \{x_0 \})\) if and only if 
\[
 u(x) = C_{k}^{H}(x-x_0) 
\]
and \(\Omega = \{ x \in \mathbb{R}^2 \mid C_{k}^{H}(x-x_0) < 1 \} \, \) for some  \(  \, k > 0 \).
\end{seur}

\textbf{Outline of this paper.}  In sections 2 and 3 we recall the definitions of an absolute minimizer, viscosity solutions of the Aronsson equation and comparison with general cones and study basic theory of general cones and functions which enjoy comparison with general cones. In section 4 we prove Theorem \ref{bhatta2}. In section 5 all lemmas needed for Theorem \ref{teoria1} are piled together and finally in section 6 Theorem \ref{teoria1} and Corollary \ref{corollaarimain} are proved.

\section{Preliminaries}
\vspace*{4mm}

In this section we recall the definitions of an absolute minimizer, viscosity solutions of the Aronsson equation and comparison with general cones. We will also introduce a new definition called general AMLE property. It is a generalization of a concept called \emph{absolute minimizing Lipschitz extension} which was introduced by Aronsson \cite{Aronsson2} to characterize infinity harmonic functions.  

In this section the Hamiltionian \( H \) is assumed to satisfy the following conditions
\begin{itemize}
 \item [(a)] \(H \in C^2(\mathbb{R}^n) \),  \( H \geq 0 \) and  \( H(0) = 0 \);
\item[(b)] \( H \) is quasi-convex: 
\[
 H(tp + (1 - t )q ) \leq \max \{ H(p), H(q)\} \qquad \text{for all}\, p,q \in \mathbb{R}^n, \,\, 0 \leq t \leq 1 ;
\]
\item[(c)] \( H\) is coersive:
\[
 \lim_{|p| \to \infty} H(p) = \infty  ;
\]
\item[(d)] 
 for any \( k \geq 0 \) the set \( \{ p \mid H(p) = k \} \) contains no interior points.
\end{itemize}

\begin{defn}
 A function \( u \in W_{loc}^{1, \infty}(\Omega)\) is an absolute minimizer of the functional 
\[
 F(v, \cdot) = || H(Dv) ||_{L^{\infty}(\cdot )} 
\]
if for every \( V \subset \subset \Omega \) and \( v \in W^{1, \infty}(V)\) with boundary values \( v = u \) on \( \partial V \), it holds
\[
 || H(Dv) ||_{L^{\infty}(V )} \geq || H(Du) ||_{L^{\infty}(V )}  .
\]
\end{defn}

\begin{defn}
An upper semicontinuous function \( u : \Omega \to \mathbb{R} \)  (abbreviated \(u \in \text{USC} (\Omega)  \))  is a viscosity subsolution of the Aronsson equation
\[
 \mathscr{A}_{H} [v](x)	= D^2 v(x) \, H_p(\nabla v (x)) \cdot H_p(\nabla v (x)) =0 
\]
if for every \( \varphi \in C^2(\Omega) \) at the local maximum point \(x_0 \) of \( u - \varphi \) we have
\[
 \mathscr{A}_{H} [\varphi](x_0) \geq 0  .
\]

Similarly a lower semicontinuous function \( u : \Omega \to \mathbb{R} \, ( u \in \text{LSC} (\Omega)) \)  is a viscosity supersolution of the Aronsson equation \(\mathscr{A}_{H} [v] = 0 \) if for every \( \varphi \in C^2(\Omega) \) at the local minimum point \(x_0\) of \( u - \varphi \) we have
\[
 \mathscr{A}_{H} [\varphi](x_0) \leq 0  .
\]
Finally, a function \(u \in C(\Omega) \) is a viscosity solution of the Aronsson equation \(\mathscr{A}_{H} [v] = 0 \) if it is both a viscosity sub- and supersolution.
 \end{defn}

\begin{defn}
For \( k > 0 \) we define a general cone \( C_{k}^{H} \) by
\[
 C_{k}^{H}(x)  \coloneqq \underset{\{H(p) = k\}}{\max} \,\, p \, \cdotp x   .
\]
\end{defn}

\textbf{Remark :} \, Since we are assuming that \( H \) is quasi-convex, the set \( \{ p \mid H(p) \leq k \}\) is convex for every \( k > 0 \). Therefore by the linear programming principle we have a very important equality
\[
 C_{k}^{H}(x) = \underset{\{H(p) \leq k\}}{\max} \,\, p \, \cdotp x   .
\]
A general cone \( C_{k}^{H} \) is therefore just a support function of the convex set \(\{p \mid H(p) \leq k\} \).

Next we list a few of the most important properties of general cones. We will use these properties frequently in future, usually without referring.

\begin{prop}
 \label{kartiolle1}
\begin{itemize}
 \item[(i)] For every \( k > 0 \), the cone \( x \mapsto C_k^H(x) \) is a convex and positively homogeneous function, i.e.,
\[
 C_k^H(\lambda x) = \lambda C_k^H(x)
\]
 for all \( \lambda > 0 \).
 \item[(ii)] \( C_k^H(x) \) is nondecreasing with respect to \( k > 0 \). Moreover \( C_k^H(x) \to C_{k_0}^H(x) \) \\ locally unifromly when \( k \to k_0 \).
\item[(iii)] \( C_k^H(0) = 0 \), \(  C_k^H \geq 0 \). Moreover, we have a triangle inequality
\[
 C_k^H(y + x) \leq  C_k^H(y) + C_k^H(x)  
\]
for all \( x,y \).

\end{itemize}
\end{prop}

\textit{Proof:} See \cite{generalcomparisonGWY}.

Let us denote \( \hat{H}(p) = H(-p) \). If \( H \) is not symmetric then \( \hat{H} \) and \( H \) are two different functions.

\begin{defn}
A function \(u \in \text{USC} (\Omega)  \) enjoys comparison with general cones from above (abbreviated \(u \in \text{CGCA} (\Omega)  \)) if for every \( V \subset \subset \Omega  , \, k > 0, \, b \in \mathbb{R} \) and \( x_0 \notin V \) from
\[
 u(x) \leq b + C_k^H(x- x_0) \qquad \text{for all} \, x \in \partial V 
\]
it follows 
\[
 u(x) \leq b + C_k^H(x- x_0) \qquad \text{for all} \, x \in V   .
\]
A function \(u \in \text{LSC} (\Omega)  \) enjoys comparison with general cones from below \( (u \in \text{CGCB} (\Omega) ) \) if for every \( V \subset \subset \Omega  , \, k > 0 , \, b \in \mathbb{R} \) and \( x_0 \notin V\) from
\[
 u(x) \geq  b - C_k^{\hat{H}}(x- x_0) \qquad \text{for all} \, x \in \partial V 
\]
it follows 
\[
 u(x) \geq b - C_k^{\hat{H}}(x- x_0) \qquad \text{for all} \, x \in V   .
\]
 Finally \( u \in C(\Omega )\) enjoys comparison with general cones  \( (u \in \text{CGC} (\Omega) ) \) if \( u \in \text{CGCA} (\Omega) \cap \text{CGCB} (\Omega)\).
\end{defn}

Next Proposition is also from \cite{generalcomparisonGWY}.

\begin{prop}
\label{kartiolle1.1}
Suppose \( u \in \text{CGCA}(\Omega) \). Then \( u \in W_{loc}^{1, \infty}(\Omega)\) and for \( x_0 \in \Omega\) and \( 0 < r < \text{dist} (x_0, \partial \Omega) \)
\[
  S_{r}^{+}(H, u, x_0) \coloneqq \inf \, \{ \, k > 0 \mid u(x) - u(x_0) \leq C_{k}^{H}(x-x_0) , \, \forall x \in \partial B_r(x_0) \, \} ,
\]
is nondecreasing with respect to \( r > 0 \). Therefore the limit
\[
 S^{+}(H, u, x_0) = \lim_{r \to 0} S_{r}^{+}(H, u, x_0)
\]
exists and \( S^{+}(H, u, x_0) = H(\nabla u(x_0))\) at the points where \( u \) is differentiable. Respectively,  if  \( u \in \text{CGCB}(\Omega) \), then \( u \in W_{loc}^{1, \infty}(\Omega)\) and 
\[
  S_{r}^{-}(H, u, x_0) \coloneqq \inf \, \{ \, k > 0 \mid  u(x) - u(x_0) \geq - C_{k}^{\hat{H}}(x-x_0) , \, \forall x \in \partial B_r(x_0) \, \} 
\]
is nondecreasing with respect to \( r > 0 \). The limit
\[
 S^{-}(H, u, x_0) = \lim_{r \to 0} S_{r}^{-}(H, u, x_0)
\]
exists and \( S^{-}(H, u, x_0) = H(\nabla u(x_0)) \) at the points where \( u \) is differentiable. 
\end{prop}

The next theorem is fundamental. The prove can be found again in \cite{generalcomparisonGWY}.

\begin{teor}
\label{abs=visco}
The following conditions are equivalent:
\begin{itemize}
\item[(i)] \(u \in C(\Omega) \) is a viscosity subsolution (supersolution) of the Aronsson equation \(\mathscr{A}_{H} [v] = 0 \),
\item[(ii)] \(u \in \text{CGCA}(\Omega) \) \, (\(u \in \text{CGCB}(\Omega) \)) .
\end{itemize}
Moreover \( u \in W_{loc}^{1, \infty}(\Omega)\) is an absolute minimizer of the functional \(F(v, \cdot) = || H(Dv) ||_{L^{\infty}(\cdot )} \) if and only if \(u \) is a viscosity solution of \(\mathscr{A}_{H} [v] = 0 \).
\end{teor}

From these three definitions, comparison with cones has turned out to be the most powerful tool in the study of solutions of Aronsson equations. Our proofs will mostly rely on this characterization. However, we will need yet another characterization.

\begin{defn}
We say that a function \( u \in W_{loc}^{1, \infty}(\Omega)\) has general AMLE (absolute minimizing Lipschitz extension) property, if for every \( V \subset \subset \Omega \) and \( \lambda \geq 0\) from condition
\[
 \underset{y,x \in \partial V}{\sup} \{ u(y) - u(x) - C_{\lambda}^H(y - x) \} \leq 0
\]
it follows that
\begin{itemize}
\item[(a)] \(
 \underset{y,x \in V}{\sup} \{ u(y) - u(x) - C_{\lambda}^H(y - x) \} \leq 0
\)

and

\item[(b)] \( || H(Du) ||_{L^{\infty}(V) } \leq \lambda \)  .
\end{itemize}

\end{defn}

\vspace*{4mm}
\textbf{Remark :} \, In fact the condition \( (b) \) follows from \( (a) \). Indeed fix a point \( x_0 \in V \) at which \( u \) is differentiable. Then for every \( e \in \partial B_1(0) \) and \( t > 0 \) the condition \( (a) \) implies 
\[
 u(x_0 + t e) - u(x_0) \leq C_{\lambda}^H(t e) = t C_{\lambda}^H( e)  .
\]
 Dividing this by \( t \) and taking a limit as \( t \to 0 \) we get 
\[
 \nabla u(x_0) \cdot e \leq C_{\lambda}^H(e) = \underset{\{H(p) \leq \lambda\}}{\max} \, p \, \cdotp e   .
\]
It follows from the convexity of the set \( \{ p \mid H(p) \leq \lambda\} \) that \(
 H(\nabla u(x_0)) \leq \lambda \, . \) This implies \((b) \).

We conclude this introductory section with the following result which is rather obvious.

\begin{teor}
 Let \( \Omega \subset \mathbb{R}^n \) be a bounded domain. The following condition are equivalent:
\begin{itemize}
\item[(i)] \( u \in \text{CGC}(\Omega)\),
\item[(ii)]  \( u \in W_{loc}^{1, \infty}(\Omega)\) has general AMLE property.
\end{itemize}
\end{teor}

\textit{Proof:} \,\( (i) \Rightarrow (ii)  \,\, \) First of all, since \( u \in CGC(\Omega)\) it is locally Lipschitz continuous by Proposition \ref{kartiolle1.1}. Let \( V \subset \subset \Omega \) and fix \( \hat{x}, \hat{y} \in V \). Suppose \( \lambda > 0 \) is such that 
\[
 \underset{y,x \in \partial V}{\sup} \{ u(y) - u(x) - C_{\lambda}^H(y - x) \} \leq 0
\]
Fix \( z \in \partial V \). Then for all \( x \in \partial V \) it holds
\begin{equation}
\label{intro1}
 u(z) -  C_{\lambda}^{\hat{H}}(x - z ) \leq u(x) \leq u(z) + C_{\lambda}^H(x - z)
\end{equation}
Since \(u \in \text{CGC} (\Omega) \) the inequality (\ref{intro1}) holds for all \( x \in V \), in particularly for \(  \hat{x} \). But the point \(z \in \partial V \) was arbitrarily chosen and therefore
\[
\underset{y,x \in \partial ( V \backslash \{ \hat{x} \})}{\sup} \{ u(y) - u(x) - C_{\lambda}^H(y - x) \} \leq 0  .
\]
Repeating the argument for a domain \(V \backslash \{ \hat{x} \} \) and for \( \hat{y} \) we have
\[
 \underset{y,x \in \partial ( V \backslash \{ \hat{x}, \hat{y}  \})}{\sup} \{ u(y) - u(x) - C_{\lambda}^H(y - x) \} \leq 0  .
\]
Especially since \(\hat{x}, \hat{y} \in \partial (V \backslash \{ \hat{x}, \hat{y}  \}) \) we have
\[
 \begin{split}
  	&u(\hat{y}) - u(\hat{x}) \leq C_{\lambda}^H( \hat{y} - \hat{x})  , \\
	&u(\hat{x}) - u(\hat{y}) \leq C_{\lambda}^H( \hat{x} - \hat{y})  .
 \end{split}
\]
Since \( \hat{x}, \hat{y} \) were arbitrary 
\[
 \underset{y,x \in  V }{\sup} \{ u(y) - u(x) - C_{\lambda}^H(y - x) \} \leq 0 . 
\]
\vspace*{4mm}

\( (ii) \Rightarrow (i)  \,\, \) Suppose that \( (i) \) doesn't hold. We may assume that  \(u \notin \text{CGCA} (\Omega) \). This means that there are \( V \subset \subset \Omega\, , \lambda > 0 , \, b \in \mathbb{R}  \) and \( x_0 \notin V \) such that 
\[
 u(x) \leq b +  C_{\lambda}^H(x - x_0)
\]
in \( \partial V \), but \( u(x) > b +  C_{\lambda}^H(x - x_0) \) for some \( x \in V \). Hence there is a domain \( W \subset V \) such that
\begin{equation}
\label{teoria0.1}
  u(x) = b +  C_{\lambda}^H(x - x_0)
\end{equation}
 for all \( x \in \partial W \) and 
\begin{equation}
\label{teoria0.2}
 u(x) > b +  C_{\lambda}^H(x - x_0)
\end{equation}
for all \( x \in W \). Equality (\ref{teoria0.1}) and triangle inequality in Proposition \ref{kartiolle1} yield 
\[
 \underset{y,x \in \partial W}{\sup} \{ u(y) - u(x) - C_{\lambda}^H(y - x) \} \leq 0  .
\]
By (ii) we have
\begin{equation}
 \label{teoria0.3}
\underset{y,x \in W}{\sup} \{ u(y) - u(x) - C_{\lambda}^H(y - x) \} \leq 0 .
\end{equation}
Fix \( x_1 \in W\). Choose \(y_1 \in  \partial W \cap [x_0,x_1]  \). This can be done since \( x_0 \notin W \). Then \( \frac{x_1 - y_1}{|x_1 - y_1|} = \frac{y_1 - x_0}{|y_1 - x_0|} \) and therefore \( C_{\lambda}^H(x_1 - y_1) +  C_{\lambda}^H(y_1 - x_0) = C_{\lambda}^H(x_1 - x_0) \) by the homogeneity of \( C_{\lambda}^H\). Using (\ref{teoria0.1}) and (\ref{teoria0.2}) we obtain
\[
  u(x_1) - u(y_1) = u(x_1) - C_{\lambda}^H(y_1 - x_0) - b > C_{\lambda}^H(x_1 - x_0)  - C_{\lambda}^H(y_1 - x_0) = C_{\lambda}^H(x_1 - y_1)
\]
which contradicts (\ref{teoria0.3}). \( \qquad \Box \)


\section{About general cones}
\vspace*{4mm}

From now on the Hamiltonian \( H \) is assumed to satisfy the following conditions:
\begin{itemize}
 \item [(H1)] \(H \in C^2(\mathbb{R}^n) \),
\item[(H2)] \( H \geq 0 \) and  \( H(0) = 0 \),
\item[(H3)] \( H \) is uniformly convex, i.e. 
\[
 \beta |\xi|^2 \geq H_{pp} \, \xi \, \cdot \xi \geq \alpha |\xi|^2  \qquad \text{for all} \,  \xi \in \mathbb{R}^n,
\]
for some \(\beta \geq  \alpha > 0 \).
\end{itemize}
Notice that since the Hamiltonian \( H \) is uniformly convex it satisfies the conditions \((b) - (d) \) in section 2. Therefore the results in Proposition \ref{kartiolle1} surely still hold. However, under the assumptions \((H1)-(H3) \) we have a lot more information on general cones. 

\begin{prop}
 \label{kartiolle2}
\begin{itemize}

\item[(i)] \( C_k^H(\cdot) \in C^2( \mathbb{R}^n \backslash \{ 0 \}) \). 
\item[(ii)] Fix \( k > 0 \). For every \( x \neq 0 \) there is unique \( p_x^k \in \{ p \mid H(p) = k \} \) such that
\[
 C_k^H(x) = p_x^k \cdot x  .
\]
The reverse spherical image map \( Y_k : \partial B_1 \to H^{-1}(k) \)
\[
 Y_k(x) = p_x^k
\]
is of class \(C^1\), one to one and onto. Furthermore for the vector \( p_x^k \) it holds 
\[
 \frac{DH(p_x^k)}{|DH(p_x^k)|} = \frac{x}{|x|}
\]
 and
\[
 D C_k^H (x) = p_x^k  .
\]
In particular, \( H(D C_k^H (x)) = k \) for all \( x \neq 0 \).
\item[(iii)] \( C_k^H(y + x) =  C_k^H(y) + C_k^H(x) \) iff \( x = \lambda y \) for \( \lambda \geq 0 \). 
\item[(iv)] Fix \( x \neq 0 \) and consider a map \( Y_x : ]0,\infty[ \to \mathbb{R}^n \),
\[
 Y_x(k) = p_x^k ,
\]
where the vector \( p_x^k  \) is defined as in part (ii). Then \( Y_x \) is a continuous path from \( ]0,\infty[ \) to \( \mathbb{R}^n \).

\end{itemize}

\end{prop}

\textit{Proof:} \, \((i) \) and \( (ii)\); Fix \( k > 0 \). Under the assumptions \((H1)-(H3) \) the set \(\{ p \mid H(p) \leq k \} \) is of class \( C^2 \) and uniformly convex. This is known to imply that \( \Sigma_k = \{ p \mid H(p) = k \} \) is a \( C^2 \) hypersurface and the spherical image map \( \, \nu : \Sigma_k \mapsto \partial B_1\)
\[
 \nu(p) = \frac{DH(p)}{|DH(p)|}
\]
is a \(C^1 \)-diffeomorphism (see \cite{schneider}, Section 2.5). Therefore \(Y_k: \partial B_1 \mapsto \Sigma_k , \, Y_k = \nu^{-1} \) is of class \(C^1\), one to one and onto.

Fix \( x \neq 0 \). Take a vector \( p_x^k \in \Sigma_k \) which gives the value
\[
 C_k^H(x) = p_x^k \, \cdot \, x .
\]
Since \( \Sigma_k \) is strictly convex the choice of this maximizing vector is unique. Moreover by the Lagrange multiplier theorem we have 
\[
 \frac{DH(p_x^k)}{|DH(p_x^k)|} = \frac{x}{|x|}  .
\]
By the definition of \( \nu \) and \(Y_k \) we have \( p_x^k = Y_k(x) \) for all \( |x|=  1 \). Using the homogeneity of the cone function we have
\begin{equation}
\label{kartiolle2kohta1}
 C_k^H(x) = Y_k \left(\frac{x}{|x|} \right) \, \cdot \, x .
\end{equation}

Since \(C_k^H \) is convex, it is locally Lipschitz continuous and therefore it is differentiable almost everywhere and \( DC_k^H \in L_{loc}^{\infty}(\mathbb{R}^n ) \). Let \( x \neq 0 \) be a point at which \(C_k^H \) is differentiable. Consider the function
\[
 \varphi(y) = C_k^H(y) - Y_k \left(\frac{x}{|x|} \right) \, \cdot \, y  .
\]
The definition of cone yields \(C_k^H(y) - p \cdot y \geq 0 \) for all \( p \in \Sigma_k \). In particular, for \( p = Y_k \left(\frac{x}{|x|} \right) \). On the other hand, by (\ref{kartiolle2kohta1}), we have \( \varphi(x) = 0 \). Therefore
\[
 D \varphi (x) = DC_k^H(x) - Y_k \left(\frac{x}{|x|} \right) = 0  .
\]
Hence 
\[
 DC_k^H(x)  = Y_k \left(\frac{x}{|x|}  \right) = p_x^k 
\]
for all \( x \neq 0 \) and \( C_k^H \in C^2(\mathbb{R}^n \backslash \{ 0\}) \).
\vspace*{4mm}

\( (iii) \) follows directly from the uniqueness of the maximizing vector \( p_x^k \).
\vspace*{4mm}

\( (iv)\) First of all the map \( Y_x(k) = p_x^k \) is well defined in \( (0, \infty) \). Let \( k_j \to k_0 > 0 \).  Then \( H(Y_x(k_j)) = H(p_x^{k_j}) = k_j \to k_0 \) and by Proposition \ref{kartiolle1}
\[
 C_{k_j}^H(x)  \to C_{k_0}^H(x)  .
\]
 Let \( (k_j) \) be a subsequence such that \(p_x^{k_j} \) converges towards some \( \hat{p} \). Then by the previous observation we have \( H(\hat{p}) = k_0 \) and 
\[
\hat{p} \cdot x = \lim_{k_j \to k_0} p_x^{k_j} \cdot x = \lim_{k_j \to k_0} C_{k_j}^H(x) = C_{k_0}^H(x)  .
\]
But since the maximizing vector is unique then \(  \hat{p} = p_x^{k_0} = Y_x(k_0) \). This proves the continuity of \( Y_x( \cdot) \). \(\qquad \Box \)


\section{Proof of Theorem 1.1.}
\vspace*{4mm}

In this section we prove Theorem \ref{bhatta2}. We do this with the help of notion called \(K\)-comparison with cones, which is a certain generalization of comparison with cones.

Here is the definition of \(K\)-comparison with cones introduced by Juutinen \cite{petri}. 

\begin{defn}
 A function \( u \in C(\Omega)\) enjoys \(K\)-comparison with cones from above if for every \( V \subset \subset \Omega \) from a condition
\[
u(x) \leq a | x - x_0| + b   \qquad \text{on} \, \, \partial V
\]
for \( x_0 \notin V \, , a \geq 0 \, , b \in \mathbb{R} \) it follows
\[
 u(x) \leq  K a | x - x_0| + b   \qquad \text{in} \, \,  V  .
\]
Similarly, a function \( u \in C(\Omega)\) enjoys \(K\)-comparison with cones from below if \( - u \) enjoys \(K\)-comparison with cones from above. Finally, we say that a function \( u \in C(\Omega)\) enjoys \(K\)-comparison with cones if it enjoys the \(K\)-comparison both from above and below. Notice that \( K \geq 1 \). 

\end{defn}

It turns out that functions which enjoy \(K\)-comparison with cones has same kind of regularity properties as infinity harmonic functions. For the proof of the next proposition see Juutinen \cite{petri}.

\begin{prop}
\label{bhattapetri}
 Let \( u \in C(\Omega)\) enjoy \(K\)-comparison with cones. Then
\begin{itemize}
\item[(i)] u satisfies the maximum and minimum principles. This means that for every \( V \subset \subset \Omega \) the maximum and the minimum value of \( u \) in \( \overline{V} \) are achieved on the boundary \( \partial V\).
 
\item[(ii)] u is locally Lipschitz. Moreover the following estimate 
\[
 |Du(x)| \leq \frac{2 K \sup_{\Omega} |u|}{\text{dist}(x,\partial \Omega)} 
\]
holds for almost every \(x \in \Omega \).

\item[(iii)] if u is non-neqative and, if \(x_0 \in \Omega \) and \( 0 < r < R < \text{dist}(x_0,\partial \Omega) \), we have
\[
 u(y) \leq u(x) e^{\frac{K |y - x|}{R-r}} \qquad \text{for all} \, \, y,x \in B_r(x_0)  .
\]
This is Harnack's inequality.  
\end{itemize}
\end{prop}
 
Next we prove an obvious result which says that a function that enjoys comparison with general cones enjoys also \(K\)-comparison with cones. This follows directly from fact that a general cone is comparable with a normal cone. We need a simple, yet important, lemma to do this.

\begin{lemm}
\label{bhatta1}
For every \( k > 0 \)  define the values 
\[
 a_k \coloneqq \min_{H(p) = k} |p| \qquad \text{and} \qquad A_k \coloneqq \max_{H(p) = k} |p|  .
\]
Under the assumptions (H1)-(H3) there exists \( K \), depending only on the Hamiltonian \(H \), such that
\[
 \frac{A_k}{a_k} \, \leq K
\]
for every \( k > 0 \).
\end{lemm}

\textit{Proof:} \, It is easy to see that from conditions \( (H1)-(H3) \) it follows
\[
 \frac{\alpha}{2} |p|^2 \leq H(p) \leq \frac{\beta}{2} |p|^2  \qquad \text{for all} \, p \in \mathbb{R}^n  .
\]
Thereby for any points \( p' , p'' \in \{ \, p \mid H(p) = k \,  \} \) we have
\[
 \frac{\alpha}{2} |p'|^2 \leq H(p') = H(p'') \leq \frac{\beta}{2} |p''|^2   .
\]
Hence the claim holds for \( K = \sqrt{ \frac{\beta}{\alpha} } \).  \( \qquad \Box \)

\begin{prop}
\label{bhattaK-vertailu}
 Let \( u \in C(\Omega)\) enjoy comparison with general cones. Then \( u \) 
enjoys \\ \(K\)-comparison with cones for some \( K = K(H) \). 
\end{prop}

\textit{Proof:} \, Let \( V \subset \subset \Omega \) and \( x_0 \notin V \, , a \geq 0 \, , b \in \mathbb{R} \) such that 
\[
 u(x) \leq a | x - x_0| + b  \qquad \text{on} \, \, \partial V  .
\]
Define 
\begin{equation}
\label{K-vertailu1}
 k_a \coloneqq \max_{|p| \leq a}  H(p)   .
\end{equation}
One sees instantly that \( B_a(0) \subset \{ p \mid H(p) \leq k_a \} \). Therefore for all \( x \in \mathbb{R}^n \)
\[
  C_{k_a}^{H}(x) = \underset{\{H(p) \leq k_a\}}{\max} \, p \, \cdotp x \geq \underset{\{ |p| \leq a \}}{\max} \, p \, \cdotp x = a |x|  .
\]
Especially \( u(x) \leq  C_{k_a}^{H}(x - x_0) + b\) on \( \partial V \). Since  \(u \in \text{CGC} (\Omega)\) then
\begin{equation}
\label{K-vertailu2}
u(x) \leq  C_{k_a}^{H}(x - x_0) + b \qquad \text{in} \, \, V  .
\end{equation}
Define next 
\[
 A \coloneqq \max_{H(p) \leq k_a} |p| 
\]
which is the same thing as \( A = \max_{H(p) = k_a} |p| \), since by \((H2) \) and \( (H3)\) the maximum is reached at the boundary. From the definition of \( A \) we conlude that \(\{ p \mid H(p) \leq k_a  \} \subset  \overline{B}_A(0) \). Therefore
\begin{equation}
\label{K-vertailu3}
  C_{k_a}^{H}(x) \leq \underset{\{ |p| \leq A \}}{\max} \,\, p \, \cdotp x = A |x| 
\end{equation}
for all \( x \in \mathbb{R}^n \). 

Now we go back to the definition of \( k_a \). From (\ref{K-vertailu1}) it is easy to see that
\[
 a = \min_{H(p) = k_a} |p| .
\]
Since \(A = \max_{H(p) = k_a} |p| \), Lemma \ref{bhatta1} yields \( A \leq K a \) for \( K = \sqrt{ \frac{\beta}{\alpha} } \).
Putting this observation together with (\ref{K-vertailu2}) and (\ref{K-vertailu3}) we get 
\[
 u(x) \leq  C_{k_a}^{H}(x - x_0) + b \leq K a |x - x_0| + b \qquad \text{in} \, \, V \, . \qquad \Box
\]

\vspace*{4mm}

We are now ready to prove the first main result. In the proof we will be using the Harnack's inequality (Proposition \ref{bhattapetri} \((iii) \) ) in a situation where \(u \) enjoys \( K \)-comparison with cones in a domain \( \Omega = B_R(x_0) \backslash  \{x_0 \} \). By using the Proposition \ref{bhattapetri} \( (iii) \)  in suitable balls, we can find a constant  \( \tilde{K} \), still depending only on \( H \), such that for all points \( x,y \in \partial B_r(x_0) \) with \( 0 < r < \frac{R}{2} \) the inequality
\[
  u(y) \leq \tilde{K} \, u(x)
\]
holds. In fact, if we do this carefully enough, we may take \(\tilde{K} = e^{K \pi} \) (see Bhattacharya \cite{Bhattacharya}).
 \vspace*{4mm}

\textbf{Proof of Theorem 1.1.} \,We may assume that \( x_0 = 0 \). For \( r \in (0,R)  \) define 
\[
 m(r) = \min_{|x| = r } u(x) \qquad \text{and} \qquad M(r) = \max_{|x| = r } u(x)  .
\]
Notice that by Proposition \ref{bhattaK-vertailu} \( u \) enjoys \(K\)-comparison with cones in \( B_R \backslash  \{0 \} \).

First we claim that there is \( r_0 \) such that both \(m(r) \) and \(M(r) \) are monotone in \( (0,r_0) \). Suppose this were untrue. Assume first that \( m(r) \) is not monotone near zero. It follows that there are radii \(0 < r_1 < r_2 < r_3 < R \) such that
\[
 m(r_1) > m(r_2) \qquad \text{and} \qquad m(r_3) > m(r_2)  .
\]
But this immediately violates the minimum principle (Proposition \ref{bhattapetri} (i)).
In case \( M(r) \) is not monotone near zero, we get a contradiction by maximum principle. Therefore the first claim holds. In particular, the limits
\[
 m_0 = \lim_{r \to 0+} m(r) \qquad \text{and} \qquad M_0 = \lim_{r \to 0+} M(r)
\]
exist (but they might not be finite).

Next we will proof that \( M_0 < \infty \). To do this, we first use Harnack's inequality (Proposition \ref{bhattapetri}(iii)) in a way discussed earlier. Thus for points \( x,y \in \partial B_r \) with \( 0 < r < R/2 \) we have the inequality
\begin{equation}
 \label{bhatta-teor-1}
u(y) \leq e^{K \pi} \, u(x) .
\end{equation}
Therefore \( M(r) \leq e^{K \pi} \, m(r) \) for all \( r \in (0, R/2) \). Hence we just need to show that \(m_0 < \infty \).

Fix \( r \in (0,R) \). For \( \rho \in (r, R ) \) choose a point \( |x_{\rho}| = \rho \) such that 
\[ 
 u(x_{\rho}) = m(\rho) ,
\]
and denote
\[
 x_r = \frac{r x_{\rho}}{\rho}  .
\]
By the definition of \( m(r) \) and non-negativity of \( u \) we have for every \( x \in \partial (B_R \backslash B_r) \)
\[
 u(x) \geq  - \left(\frac{m(r)}{R - r} \right) \, |x - x_r | + m(r)  .
\]
Since \( u \) enjoys \(K\)- comparison with cones we have 
\[
 u(x) \geq  - K \, \left(\frac{m(r)}{R - r} \right) \, |x - x_r | + m(r) 
\]
for all \( x \in B_R \backslash B_r \), especially for \( x_{\rho}\). Hence
\[
 \begin{split}
 m(\rho) = u( x_{\rho}) &\geq - K \, \left(\frac{m(r)}{R - r} \right) \, |x_{\rho} - x_r | + m(r) \\
			&= - K \, \left(\frac{m(r)}{R - r} \right) \,(\rho - r) + m(r) .
 \end{split}
\]
We write this in a slightly different way
\begin{equation}
 \label{bhatta-teor-2}
 \frac{m(\rho) - m(r)}{\rho - r} \geq - K \, \left(\frac{m(r)}{R - r} \right) 
\end{equation}
which holds for every \( 0 < r < R \) and all \( r < \rho < R \). By Proposition \ref{bhattapetri} (ii) \( u \) is locally Lipschitz in \(B_R \backslash  \{0 \} \) and therefore \( m(r)\) is also locally Lipschitz in \( (0,R) \). Therefore by taking a limit \( \rho \to r + \) in (\ref{bhatta-teor-2}) we have 
\[
 m'(r) \geq - K \, \left(\frac{m(r)}{R - r} \right)
\]
at the points where \( m(r) \) is differentiable. Dividing this by \( m(r) \) and integrating with respect to \(r\) from some small \( \tilde{r} > 0 \) to \( R/2 \) we obtain
\[
 \log \left( \frac{m(R/2)}{m(\tilde{r})} \right) \geq K \, \log  \left( \frac{R/2}{R - \tilde{r}} \right)  .
\]
Hence 
\[
 m_0 = \lim_{\tilde{r} \to 0} m(\tilde{r}) \leq \lim_{\tilde{r} \to 0} m(R/2) \, \left( \frac{R/2}{R - \tilde{r}} \right)^{-K} = 2^K \,  m(R/2) < \infty  .
\]

Now we know that the limits \( m_0 = \lim_{r \to 0+} m(r) \) and \( M_0 = \lim_{r \to 0+} M(r) \) exists and that they are finite. Next we show that 
\[
 m_0 = M_0
\]
 and we are done.  

For \( r,r' \in (0, R) \) we define
\[
 m(r,r') = \min \{ m(r), m(r') \} .
\]
Fix \( r \in (0, R/3) \) and define a function
\[
 w(x) = u(x) - m(r, 3r) .
\]
By minimum principle we have that 
\[
 w(x) \geq 0
\]
for all \( x \in B_{3r} \backslash B_r \). Since \( u \) enjoys \(K\)-comparison with cones in \( B_R \backslash \{0\}\), \( w  \)  enjoys \(K\)-comparison with cones in \( B_R \backslash \{0\}\). Hence \( w \) satisfies Harnack's inequality (\ref{bhatta-teor-1}) and therefore we have
\[
  0 \leq M(2r) - m(r, 3r) =\max_{|x|= 2r} w(x) \leq e^{K \pi} \, \min_{|x|= 2r} w(x) =  e^{K \pi} \, ( m(2r) - m(r, 3r) )  .  
\]
By sending \(r \to 0 \) and noticing that \(\lim_{r \to 0 }  m(r, 3r) = m_0 \) we get
\[
 0 	\leq M_0 - m_0 =\lim_{r \to 0 } M(2r) - m(r, 3r) \leq e^{K \pi} \, \lim_{r \to 0 } (m(2r) - m(r, 3r)) = 0 .
\]
Hence \( M_0 = m_0  . \qquad \Box \)


\section{Lemmas}
\vspace*{4mm}

In this section we list all the important Lemmas that are needed to prove Theorem \ref{teoria1}. The following fact is obvious but it is used frequently and is therefore stated separately.

\begin{lemm}
\label{lemma1} 
Let \(u \in W^{1, \infty}(\Omega) \) be such that \( || H(Du(x)) ||_{L^{\infty}(\Omega)} \, \leq k_0 \) for some \( k_0 > 0 \). Then 
\[
u(y) - u(x) \leq  C_{k_0}^{H}(y-x)
\]
for all points \( y,x \in \Omega \) for which the line segment \( [x,y] \subset \Omega \).
\end{lemm}

\textit{Proof:} \, By choosing a path \( \xi : [0,1] \rightarrow \mathbb{R} \, , \, \xi(t) = ty + (1 - t) x \) , we get

\[
\begin{split}
 u(y) - u(x) 	&= \int\limits_0^1 Du(\xi(t)) \cdot \dot{\xi}(t) \, dt \leq \int\limits_0^1 C_{\underbrace{H(Du(\xi(t)))}_{\leq k_0}}^{H}( \dot{\xi}(t)) \, dt \\
		&\leq \int\limits_0^1 C_{k_0}^{H}(y - x) \, dt =  C_{k_0}^{H}(y - x)\, . \qquad \Box
\end{split}
\]

\begin{lemm}
\label{lemma2.2}
Suppose \( u \in \text{CGCB} ( B_r(x_0) ) \) and \( || H(Du(x)) ||_{L^{\infty}( B_r(x_0))} \, \leq 1 \). If \( S^-(u,x_0) = 1 \) then there exists \( e \in \partial B_1(0)\) such that 

\[
 u(te + x_0) = u(x_0) + tC_{1}^{H}(e)
\]
for all \( t \in ]-r,0] \).
\end{lemm}

\textit{Proof:} \, We may assume that there is \( x_r \in \partial B_r(0)\) such that
\[
\begin{split}
 u(x_r + x_0) - u(x_0) 	&= -C_{S_{r}^{-}(u,x_0)}^{\hat{H}}(x_r) \\
		&\leq -C_{S^{-}(u,x_0)}^{\hat{H}}(x_r) = - C_{1}^{\hat{H}}(x_r) = - C_{1}^{H}(-x_r) ,
\end{split}
\]
where the inequality follows from Propositions \ref{kartiolle1} and \ref{kartiolle1.1}. By Lemma \ref{lemma1} for every \( 0 \leq s < 1 \) it holds
\[ 
u(s x_r + x_0) - u(x_0) \geq - C_{1}^{H}( -s x_r)  = -s C_{1}^{H}( -x_r)  . 
\] 
Therefore we have \(  u(x_r + x_0) - u(x_0) = - C_{1}^{H}(- x_r)\). On the other hand, for every \( 0 \leq s < 1 \)
\[\begin{split}
u(s x_r + x_0) - u(x_0) &= u(x_r + x_0) - u(x_0) - (\overbrace{u(x_r + x_0) - u(s x_r + x_0)}^{\geq -C_{1}^{H}(-sx_r + x_r)}) \\
		&\leq  -C_{1}^{H}(-x_r) + (1-s)C_{1}^{H}(-x_r) \\
		&= -s C_{1}^{H}(-x_r).
\end{split} \]
Thus
\[
 u(s x_r + x_0) - u(x_0) = -s C_{1}^{H}( -x_r) 
\]
holds for \( 0 \leq s < 1 \) and the claim is true for \(e = -\frac{x_r}{|x_r|} . \qquad \Box \)

\begin{lemm}
 \label{lemma3}
Let \( u \in  W^{1, \infty} ( \mathbb{R}^n ) \) and \( e \in \partial B_1(0)\) be such that
\begin{itemize}
\item[(i)] \( || H(Du) ||_{L^{\infty}(\mathbb{R}^n ) } \leq 1 \)
\item[(ii)] \( u(te) = t C_1^{H}(e) = t p_e \cdot e \) \, for all \( t \geq 0 \)
\end{itemize}
where \( p_e \in  H^{-1}( 1 )  \) is such that \( C_1^{H}(e) = p_e \cdot e \). Then 
\[ 
u(x) \geq p_e \cdot x 
 \]
 for all \( x \in \mathbb{R}^n \). Moreover, if in part \( (ii) \) we have an equality for all \(t \in \mathbb{R}\) then \( u(x) = p_e \cdot x \) for all \( x \).
\end{lemm}

\textit{Proof:} \, Fix \( x \in \mathbb{R}^n  \). Assumptions (ii) and (i) yield 
\begin{equation}
\label{lemma3.1} 
 u(x) - t C_1^{H}(e) = u(x) - u( t e) \geq - C_1^{H}(t e - x)
\end{equation}
for all \( t \geq 0 \). For \( t > 0 \) let \( p_{e - \frac{x}{t}} \in H^{-1}( 1 ) \) be a vector such that 
\[
 C_1^{H}( e - \frac{x}{t}) = p_{e - \frac{x}{t}} \cdot (e - \frac{x}{t})  .
\]
Then by (\ref{lemma3.1} ) for all \( t \geq 0 \) 
\[
\begin{split}
u(x) 	&\geq t C_1^{H}(e) - C_1^{H}(t e - x) = tC_1^{H}(e) - t C_1^{H}( e - \frac{x}{t}) \\
	&= t C_1^{H}(e) - t \, p_{e - \frac{x}{t}} \cdot ( e - \frac{x}{t} ) = t \underbrace{( \max_{H(p)=1} p \cdot e - p_{e - \frac{x}{t}} \cdot e)}_{\geq 0} + p_{e -\frac{x}{t}}\cdot x \\
	&\geq p_{e -\frac{x}{t}} \cdot x  .
\end{split}
\]
Therefore \(  u(x) \geq \lim_{t \to \infty} p_{e -\frac{x}{t}} \cdot x = p_e \cdot x \). The other part of the lemma follows from repeating the argument for a function \( \tilde{u}(x) = -u(-x)  \) .  \( \qquad \Box \)

Next Lemma is crucial in the proof of Theorem \ref{teoria1}.

\begin{lemm}
 \label{lemma4}
Let \( n \geq 2 \). Suppose \( u \in  W^{1, \infty} ( \mathbb{R}^n ) \) is such that
\begin{itemize}
\item[(i)] \( || H(Du) ||_{L^{\infty}(\mathbb{R}^n ) } \leq 1 \),
\item[(ii)] \( u(0) = 0 \) and there exists \(\epsilon > 0 \) such that \( u(x) \geq  - C_{1 - \epsilon}^{\hat{H}}(x)\) for all \( x \in \mathbb{R}^n  \),
\item[(iii)] \(u \) is a viscosity supersolution of the Aronsson equation in \(\mathbb{R}^n \backslash \{ 0 \} \)
\item[(iv)] there exists \( \hat{e} \in \partial B_1(0) \) such that
\[
  u(t \hat{e} ) = t C_{1}^{H}( \hat{e} )
\]
for all \( t \geq 0 \). 

Then 
\[
 u(x) = C_{1}^{H}(x) .
\]

\end{itemize}

\end{lemm}

\textit{Proof:} \, First define a set of functions
\[
 \mathscr{K} \coloneqq \{ v \in W^{1, \infty} ( \mathbb{R}^n ) \mid v \, \text{satisfies conditions} \, (i) - (iv) \,\text{with uniform} \, \epsilon > 0 \, \text{in part} \, (ii) \}
\]
and a function \( w(x)\coloneqq \underset{ v \in \mathscr{K} }{\inf} v(x) \) . By Ascoli's theorem \( w \in \mathscr{K} \). Since for \(\lambda > 0, \,  \frac{w( \lambda x)}{\lambda} \in \mathscr{K} \), and thus
\[ 
 w( x) \leq \frac{w( \lambda x )}{\lambda}  . 
\] 
It follows that \( w \) is positively homogeneous meaning
\[
\label{homogeeninen}
w(  \lambda x) = \lambda w (x)
\]
for every \( \lambda \geq 0 \). Furthermore, 
\begin{equation}
\label{pienempi kuin kartio}
 w(x) \leq u(x) \leq C_{1}^{H}(x),
\end{equation}
where the second inequality follows from assumption (i). The plan of the proof is to look at the set
\[
 F = \{ x \in \mathbb{R}^n \backslash \{ 0 \} \mid w(x) = C_{1}^{H}(x) \}  .
\]
Surely \(F\) is closed. The goal is to show that \( F\) is also open and therefore by assumption \( (iv)\) it has to be the whole \( \mathbb{R}^n \backslash \{ 0 \} \). The proof is then completed by (\ref{pienempi kuin kartio}). 

Fix \( x_0 \in F\). Since both \( w \) and the cone \( C_{1}^{H} \) are positively homogeneous functions, we have
\[ 
w ( t e ) = t C_{1}^{H}(e)
\]
for all \( t \geq 0 \), where \( e = \frac{x_0}{|x_0|}\). Thereby assumption \( (i) \) and Lemma \ref{lemma3} yield
\begin{equation}
 \label{alhaalta rajoitettu}
w(x) \geq p_e \cdot x
\end{equation}
for all \( x \in  \mathbb{R}^n \) where the vector \(p_e \) is such that \(C_{1}^{H}(e) = p_e \cdot e \). Denote \( T \coloneqq \{ x \in  \mathbb{R}^n \mid p_e \cdot x = 0  \}\) .

\begin{figure}[t]
\centering
\psfrag{e}{$e$}
\psfrag{x}{$- \bar{x}$}
\psfrag{z0}{$z_0$}
\psfrag{t0e}{$t_0 e$}
\psfrag{Br}{$B_r(z_0)$}
\includegraphics[scale=0.65]{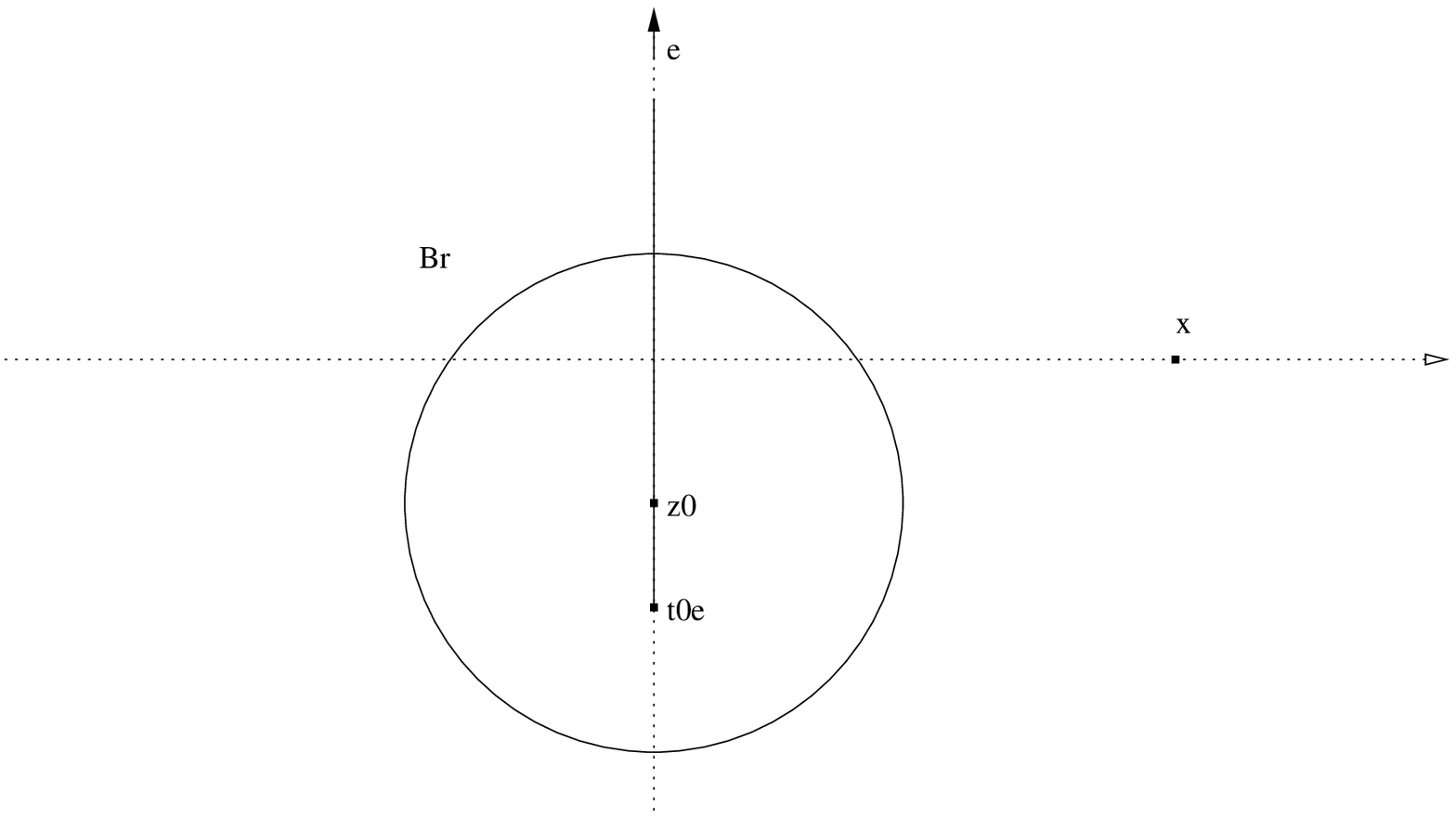}
\caption{}
\label{kuva1}
\end{figure}

\textbf{Step 1 :} \, We show that 
\begin{equation}
 \label{apuvaite}
\underset{ \{ x \in T, \, |x|=1 \}}{\min} w(x) > 0  .
\end{equation}
Suppose this were not true and there would be \(  \bar{x} \in T,  |\bar{x}|=1 \) such that \( w(\bar{x}) = 0 \). Then the function 
\[  
h(x) = w(x +\bar{x}) 
\] 
is a supersolution of \(  \mathscr{A}_{H}[v]= 0 \) in \(\mathbb{R}^n \backslash \{-\bar{x} \}\), satisfies the condition \( (i) \) and by (\ref{alhaalta rajoitettu}) we have
\[
 h(x) \geq p_e \cdot x  .
\]
Assumption (i) yields
\[
 h(t e) = h(te) - h(0) \leq C_{1}^{H}(t e) = t p_e \cdot e  .
\]
for all \( t \geq 0 \). Hence 
\[  
h(t e) = t p_e \cdot e = t C_{1}^{H}( e) 
\] 
for all \( t \geq 0 \). Denote
\[
 t_0 = \inf \, \{ t \in \mathbb{R} \mid h(s e) = s C_{1}^{H}( e) \, \, \text{for all}  \, \, s \geq t \}  .
\]
By the previous discussion \( t_0 \leq 0 \). Suppose \( t_0 > - \infty \). Then there would be \( r > 0 \) such that \( B_{2r}(t_0e) \subset \mathbb{R}^n \backslash \{-\bar{x} \}\). For \( z_0 = (t_0 + \frac{r}{2})e \) the assumption \((iii) \) yields \( h \in \text{CGCB}(B_{r}(z_0)) \). Moreover we have \( || H(Dh(x)) ||_{L^\infty} \leq 1 \) and for all \( - \frac{r}{2} < t < r \) it holds
\begin{equation}
\label{apu1}
 h( te + z_0) = h(z_0) + t C_{1}^{H}(e)
\end{equation}
as drawn in Figure \ref{kuva1}. This implies \( S^-(h,z_0) = 1 \). Then by Lemma~\ref{lemma2.2} 
\begin{equation}
\label{apu2}
 h( t \tilde{e} + z_0) = h(z_0) + t C_{1}^{H}(\tilde{e})
\end{equation}
for all \( t \in ]-r, 0[\) for some \( |\tilde{e}| = 1\). From (\ref{apu1}) and (\ref{apu2}) it follows
\[
 h(re + z_0) - h(-r \tilde{e} + z_0) = r C_{1}^{H}(e)+ r C_{1}^{H}(\tilde{e}). 
\]
On the other hand \(h(re + z_0) - h(-r \tilde{e} + z_0) \leq  C_{1}^{H}( re+r \tilde{e})\). Therefore we have
\[
 C_{1}^{H}(e)+ C_{1}^{H}(\tilde{e}) = C_{1}^{H}( e+ \tilde{e})
\]
which can only be true if \(\tilde{e} = e \) by Proposition \ref{kartiolle2}. Hence (\ref{apu2}), together with the definitions of \( t_0 \) and \(z_0 \), yields
\[
h( te ) = t  C_{1}^{H}(e) 
\]
for all \( t \geq t_0 - \frac{r}{2} \) which is a contradiction. Therefore \( t_0 = - \infty\) and 
\[
 h(te) =  t C_{1}^{H}(e)
\]
for all \( t \in \mathbb{R} \). Thereby Lemma \ref{lemma3} yields 
\[
 h(x) = p_e \cdot x .
\]
We also have \( w(x) = p_e \cdot x  \) (remember that \(h(x) = w(x +\bar{x})\) and \( p_e \cdot \bar{x} = 0\) ) . But this contradicts the assumption \( (ii) \). Therefore (\ref{apuvaite}) must hold.   

\vspace*{4mm}
Using the observation (\ref{apuvaite}) and the homogeneity of \(w \) we find \( \delta > 0 \) such that 
\begin{equation}
\label{gammanmaar}
w > 0 \, \, \text{in} \, \, \Gamma_{\delta}
\end{equation}
where 
\[ 
\Gamma_{\delta} = \{ x \in \mathbb{R}^n \mid p \cdot x > 0 \, \, \text{for some}\, p \, \text{for which}\, |p - p_e|< \delta \}  .
 \]

\vspace*{6mm}
\textbf{Step 2 :} \, We define a quantity \( k(x) \) to be a positive number for which
\[
  C_{k(x)}^{H}(x) = w(x)  .
\]
Since \( w > 0 \) in \( \Gamma_{\delta} \) by (\ref{gammanmaar}) and \( k(x) \) is uniquely determined, the map  \(x \mapsto k(x) \) is well defined in \( \Gamma_{\delta} \). Moreover, since \( w \) is continuous, the map \(x \mapsto k(x) \) has to be continuous as well. We also find a unique vector \( p_x^{k(x)}\) on the set \( H^{-1}(k(x)) \) which gives the value 
\[ 
 C_{k(x)}^{H}(x) =  p_x^{k(x)} \cdot x . 
\]
Then by Lemma \ref{kartiolle2} the map 
\[
x \mapsto p_x^{k(x)} 
 \]
is also continuous in \( \Gamma_{\delta} \). Since \( w(e) = C_{1}^{H}(e) \), we have \( k(e) = 1\). In particular, \( p_x^{k(x)} \rightarrow p_e \) as \( x \rightarrow e \).

\begin{figure}[t]
\centering
\psfrag{pe}{$p_e$}
\psfrag{px}{$- p_x^{k(x)}$}
\psfrag{G}{$\Gamma_{\delta}$}
\psfrag{V}{$V_x$}
\psfrag{x}{$x$}
\includegraphics[scale=0.65]{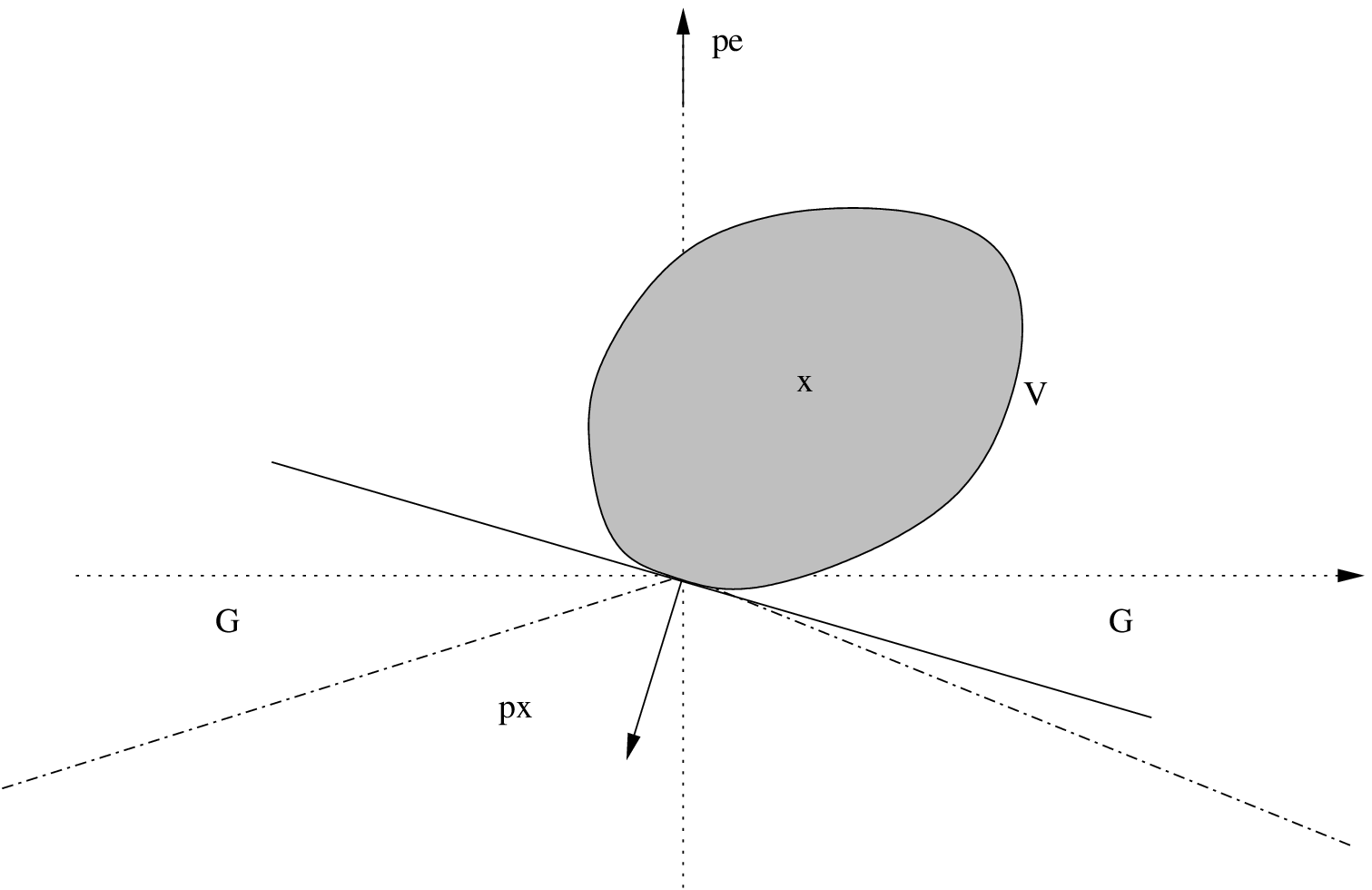}
\caption{}
\label{kuva2}
\end{figure}

Fix \( x \in \Gamma_{\delta} \). Define a set
\[
 V_x \coloneqq \{ y \in\mathbb{R}^n \mid   C_{k(x)}^{\hat{H}}(y - x) <  w(x)  \}  .
\]
From the definition of \( k(x) \) it follows directly that \( 0 \in \partial V_x \) for every \( x \). The boundary \( \partial V_x \) is a smooth hypersurface and its outer normal at the point \( 0 \in \partial V_x \) is \( D C_{k(x)}^{\hat{H}}(- x) \). Moreover, we have
\[ 
 D C_{k(x)}^{\hat{H}}(- x) = - D C_{k(x)}^{H}( x) = - p_x^{k(x)}   .
\]
It follows from the convexity of \( V_x \) that the tangent space \(\{y \in \mathbb{R}^n \mid (- p_x^{k(x)}) \cdot y = 0  \} \) supports \( V_x \) at the origin. Hence
\begin{equation}
 \label{inkluusio1}
 V_x 	\subset \{ y \in \mathbb{R}^n \mid  p_x^{k(x)} \cdot y > 0 \}  .
\end{equation}
The picture is as shown in Figure \ref{kuva2}.

Since  \( p_x^{k(x)} \rightarrow p_e \) as \( x \rightarrow e \) there is some small neighbourhood \( U_e \) of \(e \) such that \( | p_x^{k(x)} - p_e| < \delta \) for all \( x \in U_e \), where \( \delta \) is the same as in  (\ref{gammanmaar}). This together with (\ref{inkluusio1}) and the definition of \( \Gamma_{\delta}  \) implies
\[ 
 V_x \subset \Gamma_{\delta} 
\]
for all \(x \in U_e  \). By (\ref{gammanmaar})
\begin{equation}
\label{w_positiivinen}
w > 0 \, \, \text{in} \, \, V_x
\end{equation}
for every \( x \in U_e \).

Fix \( x \in U_e \). From the definition of \(V_x\) it follows that \(  w(x) - C_{k(x)}^{\hat{H}}(y - x) = 0 \) for all \( y \in \partial V_x  \). On the other hand, by (\ref{w_positiivinen}), we have \( w(y) \geq 0 \) for every \( y \in \partial V_x  \). Since \( u \in \text{CGCB}(\mathbb{R}^n \backslash \{ 0 \} ) \) we have
\[ 
w(y) \geq w(x) - C_{k(x)}^{\hat{H}}(y - x)
\]
for all \(y \in V_x \). In particular, when \(r > 0 \) is so small that \( B_r(x) \subset V_x \), then 
\[
 w(y) - w(x) \geq  - C_{k(x)}^{\hat{H}}(y - x)
\]
for all \( y \in \partial B_r(x) \). By the defintion of \(S_r^-(w, x) \) we have
\[
 k(x) \geq S_r^-(w, x) \, .
\]
Proposition \ref{kartiolle1.1} yields 
\[
   k(x) \geq S^-(w, x) = H(Dw(x)) 
\]
for all points \(x \in U_e \) where \( w \) is differentiable. Combining this with the definition of \( k(x) \) we have 
\begin{equation}
\label{epayhtalo1}
 w(x) = C_{k(x)}^{H}(x) \geq C_{H(Dw(x))}^H(x)
\end{equation}
for almost every \(x \) in \( U_e \). 

On the other hand, it follows from the homogeneity of \( w \) that
\begin{equation}
\label{epayhtalo2}
w(x) = Dw(x) \cdot x \leq C_{H(Dw(x))}^H(x)  .
\end{equation}
Combining (\ref{epayhtalo1}) and (\ref{epayhtalo2}) we have that 
\begin{equation}
\label{askel2.1}
  w(x) = C_{H(Dw(x))}^H(x)  ,
\end{equation}
\( k(x) = H(Dw(x) ) \) and 
\begin{equation}
\label{askel2.2}
 Dw(x) = p_x^{k(x)} 
\end{equation}
for almost every \( x \in U_e \).

\vspace*{6mm}
\textbf{Step 3 :} \, The last thing we need to show is that \( H(Dw(x) ) \equiv 1 \). By the convexity of \( H\) we have
\begin{equation}
\label{askel3.1}
  DH (p) \cdot (q - p) \leq H(q) - H(p) \leq  DH (q) \cdot (q - p) \, .
\end{equation}
for all \( p, q \in \mathbb{R}^n\). We will use this later.

From (\ref{askel2.2}) we conclude two things. First, \( w \in C^1(U_e) \). This follows from the continuity of the map \(x \mapsto p_x^{k(x)} \). Second, by Lemma \ref{kartiolle2} we have the relation
\[
\frac{DH(Dw(x))}{|DH(Dw(x))|} = \frac{x}{|x|} 
\]
in \( U_e \). Fix now \( x \in U_e \) and \( |\eta | = 1 \). Using these observations and the homogeneity of \(w\) ( the fact that \( w(x) = Dw(x) \cdot x  \) ) we have for \(  \lambda = \frac{|DH(p_x^{k(x)})|}{|x|} \)
\[
\begin{split}
 \lim_{t \to 0 }  DH(Dw(x)) &\cdot \left( \frac{Dw(x + t \eta) - Dw(x)}{t} \right) \\
	&= \lambda \lim_{t \to 0 } \frac{Dw(x + t\eta) \cdot x - Dw(x) \cdot x }{t} \\
	&= \lambda \left( \lim_{t \to 0 } \frac{Dw(x + t\eta) \cdot (x + t \eta) - Dw(x) \cdot x }{t} - Dw(x + t \eta) \cdot \eta  \right) \\
	&= \lambda \left( \lim_{t \to 0 }  \frac{w(x + t\eta) - w(x)}{t} - Dw(x + t \eta) \cdot \eta   \right)  \\
	&= 0 .
\end{split}
\]
Similarly we conlude that 
\[
 \lim_{t \to 0} DH(Dw(x + t \eta)) \cdot \left( \frac{Dw(x + t \eta) - Dw(x)}{t} \right) = 0  . 
\]
Therefore by using (\ref{askel3.1}) for \( p = Dw(x) \) and \(q = Dw(x + t \eta) \) we obtain
\[
 \lim_{t \to 0} \frac{H(Dw(x + t \eta) ) - H(Dw(x))}{t} = 0
\]
for every \( |\eta| = 1 \). This means that the function \( H(Dw( \cdot)) \) is differentiable and its gradient vanishes at every point \(x \in U_e\). Therefore 
\[
 H(Dw(x)) \equiv \, \text{constant} \, = 1
\]
in \( U_e \). By (\ref{askel2.1}) we have
\[
 w(x) = C_1^H(x)
\]
in the neighbourhood \( U_e \) of \(e \).

However, our goal was to show that the chosen point \( x_0 \), not \( e \), is an interior point of
\[
  F = \{ x \in \mathbb{R}^n \backslash \{ 0 \} \mid w(x) = C_{1}^{H}(x) \}  .
\]
But this follows from the fact that \( e = \frac{x_0}{|x_0|}\) and \(w \) and the cone \( C_{1}^{H} \) are homogeneous functions. Hence \( x_0 \) is an interior point of \( F \) which is therefore an open set. \( \qquad \Box\)

\vspace*{4mm}

The last lemma is a general result for viscosity solutions. We sketch it only for Aronsson equations of type (\ref{aronsson2}), although a similar result is true for more general elliptic equations. For further details see \cite{SaWaYu} (Lemma 4.2).

\begin{lemm}
\label{viimeinenlemma}
Let \( u \in C(B_1) \) be a viscosity solution of 
\[
 \mathscr{A}_{H} [u] = 0 \qquad \text{in} \, \, B_1 \backslash \{ 0 \}  .
\]
Then one of the following holds:
\begin{itemize}
\item[(i)] \( u \) is a viscosity solution in the whole ball;
\item[(ii)] there exists \( \epsilon > 0 \) and \( p \neq 0 \) such that  
\[
 u(x) \geq u(0) + p \cdot x + \epsilon | x | \qquad \text{in} \, B_{\epsilon}  ;
\]
\item[(iii)] there exists \( \epsilon > 0 \) and \( q \neq 0 \) such that
\[
 u(x) \leq u(0) + q \cdot x - \epsilon | x | \qquad \text{in} \, B_{\epsilon}  .
\]
\end{itemize}
\end{lemm}

\textbf{Remark:} \, It is easy to see that in case \( (ii) \) \( u \) is a viscosity subsolution but not a supersolution in the whole ball. On the other hand, in case  \( (iii) \) \( u \) is a viscosity supersolution but not a subsolution in \(B_1 \).


\section{Proof of Theorem 1.2 and Corollary 1.3}
\vspace*{4mm}

\textbf{Proof of Theorem 1.2. } \, We may assume that \( x_0 = 0 \). By Theorem \ref{bhatta2} \( u \) can be extended continuously to the whole ball. We denote
\[
 b = u(0) = \lim_{x \to 0} u(x)  . 
\]
From now on we may assume that \(b = 0 \). Suppose that \( (i) \) doesn't happen. By Lemma \ref{viimeinenlemma} \(u \) is either viscosity super- or subsolution of the Aronsson equation in the whole ball. Suppose it is a subsolution, but not a supersolution. Then by Lemma \ref{viimeinenlemma} there exists \( p_0 \neq 0 \) and \( \epsilon > 0 \) such that
\[
u(x) \geq p_0 \cdot x + \epsilon |x| \qquad \text{in} \, B_{\epsilon}  .
 \]
Because of this there is \( \delta > 0 \) and a smaller neighbourhood of origin \( V \subset \subset B_r(0) \) such that
\begin{equation}
 \label{teoria1.2}
u(x) \leq p_0 \cdot x + \delta
\end{equation}
for all \( x \in  V \) and 
\begin{equation}
 \label{teoria1.3}
u(x) = p_0 \cdot x + \delta
\end{equation}
 for all \( x \in \partial V \).

Denote \( \lambda \coloneqq || H(Du(x))||_{L^{\infty}(V)} \). Choose \( e = \frac{DH(p_0)}{|DH(p_0)|}\). By Proposition \ref{kartiolle2} we have
\[
 C_{H(p_0)}^H(e)  = p_0 \cdot e   .
\]
Denote 
\[
 \bar{t} = \sup \, \{ t > 0 \mid [0, te] \subset V \}  .
\]
Then \( \bar{t}e \in \partial V  \) and therefore by Lemma \ref{lemma1} and (\ref{teoria1.3}) we have
\[
  C_{\lambda}^H(\bar{t}e) \geq u( \bar{t}e) - u(0) = \bar{t} \, p_0 \cdot e + \delta  = \bar{t} \,  C_{H(p_0)}^H(e) + \delta > C_{H(p_0)}^H( \bar{t} e) .
\]
Therefore 
\begin{equation}
\label{teoria.lisays}
 \lambda > H(p_0).
\end{equation} 
Furthermore, for all \( x \in V \) it holds
\begin{equation}
  \label{teoria1.4}
u(x) 	\geq p_0 \cdot x + \epsilon |x|\geq - C_{H(p_0)}^{\hat{H}}( x) + \epsilon |x| \geq - C_{ (\lambda - \epsilon ' ) }^{\hat{H}}( x)
\end{equation}
for some \( \epsilon ' > 0 \). Consider next the quantity
\[
 \underset{x \in \partial V}{\max} \{ u(x) -C_{\lambda}^H(x) \}  .
\]
Suppose this maximum were strictly smaller that zero. Then there would be \( \lambda ' \in ]H(p_0), \lambda[ \) such that
 \[
 \underset{x \in \partial V}{\max} \{ u(x) -C_{\lambda '}^H(x) \} \, \leq 0 .
\] 
On the other hand, by (\ref{teoria1.3}), for all \( y,x \in \partial V \) it holds 
\[
 u(y) - u(x) = p_0 \cdot (y - x) \leq C_{H(p_0)}^H(y - x ) < C_{\lambda '}^H(y - x )  .
\]
Similarly for all \( x \in \partial V \) we have \(  u(0) - u(x) = - p_0 \cdot x - \delta < C_{\lambda '}^H( - x ) \, .\)  Hence 
\[
 \underset{y,x \in \partial ( V \backslash \{ 0\}) }{\max} \{u(y) - u(x) -C_{\lambda'}^H(y -x) \} \,\leq 0   .
\]
By the general AMLE property of \( u \) we would have \( ||H(Du)||_{L^{\infty}(V)} \leq \lambda ' < \lambda \) which contradicts the definition of \( \lambda \). Hence
\begin{equation}
 \label{teoria1.45}
 u(\bar{x}) - C_{\lambda}^H(\bar{x}) = \underset{x \in \partial V}{\max} \{ u(x) -C_{\lambda}^H(x) \} \, \geq 0
\end{equation}
for some \(\bar{x} \in \partial V \).

Next we claim that \( [0,\bar{x}] \subset \overline{V} \). Indeed, if this weren't true there would be \( \bar{y} \in [ 0 , \bar{x}[ \cap \partial V \) such that \( [0,\bar{y}]  \subset \bar{V} \). Then by Lemma \ref{lemma1} \( u(\bar{y}) \leq C_{\lambda}^H(\bar{y}) \). But then again by (\ref{teoria1.3}) and (\ref{teoria.lisays}) 
\[
 \begin{split}
  u(\bar{y}) 	&= u(\bar{x}) - ( u(\bar{x}) - u(\bar{y}) ) \geq  C_{\lambda}^H(\bar{x}) - p_0 \cdot ( \bar{x} - \bar{y}) \\
		&\geq  C_{\lambda}^H(\bar{x}) -  C_{H(p_0)}^H(\bar{x} - \bar{y} ) >  C_{\lambda}^H(\bar{x}) -  C_{\lambda}^H(\bar{x} - \bar{y} ) =  C_{\lambda}^H(\bar{y})
 \end{split}
\]
which is a contradiction. 

Since \( [0,\bar{x}] \subset \overline{V} \) and \( \lambda = ||H(Du)||_{L^{\infty}(V)} \) we have 
\[ 
u (t \bar{x} ) \leq C_{\lambda}^H( t \bar{x})
\] 
for all \( t \in [0,1] \).  By (\ref{teoria1.45}), \(u ( \bar{x} ) = C_{\lambda}^H(  \bar{x})\). Therefore for all \( t \in [0,1] \) we have 
\[
 u (t \bar{x} ) = u ( \bar{x} ) - (\overbrace{u(\bar{x}) - u(t \bar{x})}^{\leq C_{\lambda}^{H}(\bar{x}- t \bar{x})}) \geq C_{\lambda}^{H}(\bar{x}) - ( 1-t) C_{\lambda}^{H}(\bar{x}) = t C_{\lambda}^{H}(\bar{x}) .
\]
Hence
\begin{equation}
 \label{teoria1.5}
 u (t \bar{x} ) = C_{\lambda}^H( t \bar{x})
\end{equation}
for all \( t \in [0,1] \). 

Consider now a sequence \( (h_k) \) which has the property that \( h_k \rightarrow 0 \) as \( k \to \infty \) and 
\[
 \lim_{k \to \infty} \frac{u(h_k x)}{h_k} = w(x)
\]
locally uniformly in \( \mathbb{R}^n \). By (\ref{teoria1.4}) and (\ref{teoria1.5}) one sees that \( w \) satisfies the conditions:

\begin{itemize}
\item[(i)] \( || H(Dw) ||_{L^{\infty}(\mathbb{R}^n ) } \leq \lambda \),
\item[(ii)] \( w(0) = 0 \) and \( w(x) \geq  - C_{\lambda - \epsilon}^{\hat{H}}(x)\) for all \( x \in \mathbb{R}^n  \),
\item[(iii)] \(w \) is a viscosity solution of the Aronsson equation in \(\mathbb{R}^n \backslash \{ 0 \} \),
\item[(iv)]  for \( \hat{e} = \frac{\bar{x}}{|\bar{x}|}  \) it holds
\[
  w(t \hat{e} ) = t C_{\lambda}^{H}( \hat{e} )
\]
for all \( t \geq 0 \). 
\end{itemize}
Function \( w \) satisfies the assumptions of Lemma \ref{lemma4} with \( H \) replaced by \( H_{\lambda} = \frac{H}{\lambda} \). Lemma \ref{lemma4} yields
\[
 w(x) 	= C_1^{H_{\lambda}}( x) = C_{\lambda}^{H}(x)  .
\]
Hence the limit doesn't depend on the chosen sequence \( (h_k) \) and therefore
\[
 \lim_{h \to 0} \frac{u(h x)}{h} = C_{\lambda}^{H}( x)  .
\]
 This implies that \( u(x) = C_{\lambda}^{H}(x)  + o(|x|)  .  \qquad \Box\)

\vspace*{4mm}

 \textbf{Remark :} \, In the beginning of the proof we assumed that \( u \) is not a supersolution in the ball. If instead we assume that \( u \) is not a subsolution in the ball, we have
\[
  u(x) = u(x_0) - C_{\lambda}^{\hat{H}}(x-x_0)  + o(|x - x_0|) 
\]
for some \( \lambda > 0 \).

\vspace*{6mm}

Let us turn our attention to the last main result which is Corollary \ref{corollaarimain}. In the proof we use the result by C. Wang and Y. Yu \cite{WaYu} who showed that solutions of the Aronsson equation under the assumptions (H1) - (H3) are \( C^1 \)  (see also \cite{savin}). Moreover we have the following uniform estimate. 

\begin{teor}[\cite{WaYu}]
\label{wangyu1}
 Let \( n = 2 \). Suppose \( u \) is a solution of the Aronsson equation \( \mathscr{A}_H (u) = 0 \) in \(B_1(0) \subset \mathbb{R}^2 \), then \( u \) is differentiable at \( 0 \) and for every \( \epsilon > 0 \) there exists \( \delta = \delta( \epsilon, H) > 0 \) such that whenever
\[
  || u(x) - u(0) - e_1 \cdot x ||_{L^{\infty}(B_1)} \leq \delta,
\]
where \( 1 \leq H(e_1) \leq 2 \), then
\[
 | D u (0) - e_1 | \leq \epsilon  .
\]
\end{teor}

Using this we can derive following result.

\begin{teor}
\label{wangyu2}
 Let \( n = 2 \). For every \( \epsilon > 0 \) there is \( \delta = \delta( \epsilon, H) > 0 \) such that if we have two solutions \( u , v \) of the Aronsson equation \( \mathscr{A}_H (u) = 0 \) in \(B_1(0) \subset \mathbb{R}^2 \) for which \( ||u ||_{L^{\infty}(B_1)}, ||v ||_{L^{\infty}(B_1)} \leq C \) and
\[
 || u - v ||_{L^{\infty}(B_1)} \leq \delta , 
\]
then 
\[
  | D u (0) - D v (0) | \leq \epsilon  . 
\]
\end{teor}
The proof is pretty standard, but we sketch it here for the reader's convenience. 
\vspace*{4mm}

\textit{Proof:} Suppose this would not be true. Then there would be \( \epsilon_0 > 0 \) and sequences \( u_j \) and \(v_j \) of solutions of  \( \mathscr{A}_H (u) = 0 \) in \(B_1(0) \) for which 
\begin{itemize}
 \item[(i)] \( ||u_j ||_{L^{\infty}(B_1)}, ||v_j ||_{L^{\infty}(B_1)} \leq C \) 
 \item[(ii)]  \( || u_j - v_j ||_{L^{\infty}(B_1)} \to 0 \)
\item[(iii)] \( | D u_j (0) - D v_j (0) | \geq \epsilon_0 > 0  . \)
\end{itemize}

By Proposition \ref{bhattapetri} (ii)
\[
 |Du_j(x)| \leq 4K \, ||u_j||_{L^{\infty}(B_1)} \leq 4K \, C 
\]
for almost every \( x \in B_{1/2}\). Hence \( (u_j) \) is equicontinuous in \(B_{1/2}\). We may assume \( u_j \to u \) uniformly in \( B_{1/2} \) and \( u \) is a solution of   \( \mathscr{A}_H (u) = 0 \) in \(B_{1/2} \). It follows from Theorem \ref{wangyu1} that \( u\) is differentiable at \( 0 \).   
\vspace*{3mm}

\textit{Case 1:} \(D u(0) = 0 \). Let \( \epsilon > 0 \). Then there is \(r > 0 \) such that
\[
 ||u(x) - u(0) ||_{L^{\infty}(B_r)} \leq \epsilon r 
\]
by the differentiability. When \(j \) is large enough
\[
 ||u_j(x) - u_j(0)||_{L^{\infty}(B_r)} \leq 2 \epsilon r 
\]
by the uniform convergence. Apply Proposition \ref{bhattapetri} (ii)
to the function \( w_j(x) = u_j(x) - u_j(0) \) to obtain
\begin{equation}
\label{teoria3.0}
  |Du_j(0)| 	\leq  ||Dw_j ||_{L^{\infty}(B_{r/2})} \leq \frac{4K}{r} \, || w_j ||_{L^{\infty}(B_r)} \leq 8K \, \epsilon 
\end{equation}
when \(j \) is large. Hence \( D u_j(0) \to 0 \). Using the assumption (ii) and the argument (\ref{teoria3.0}) for \( w_j(x) = v_j(x) - v_j(0)\) we have that \(  D v_j(0) \to 0 \). But this contradicts \( (iii)\). 
 
\vspace*{3mm}

\textit{Case 2:} \(D u(0) = e \neq 0 \). First we change \( H \) to \( H_{\alpha} = \alpha H \) where \( \alpha = \frac{1}{H(e)} > 0 \). By doing this we get the desired property
\[
 H_{\alpha}(e) =1 .
\]
Our new \(H_{\alpha} \) still satisfies the conditions (H1) - (H3) and if \(w \) is any solution of \( \mathscr{A}_H (w) = 0 \) then it is also a solution of 
\[ 
 \mathscr{A}_{H_{\alpha}} (w) = 0 .
\] 
Fix now \( \delta = \delta (\epsilon_0 /3, H_{\alpha}) \) in Theorem \ref{wangyu1}. Denote 
\( u_r(x) = \frac{u(rx)- u(0)}{r}  ,  \, u_{j,r}(x) = \frac{u_j(rx) - u_j(0)}{r} \) and \(v_{j,r}(x) = \frac{v_j(rx) - v_j(0)}{r} \). Using differentiability of \( u\) we find \( r > 0 \) such that
\begin{equation}
\label{teoria3.1}
  ||u_r(x) - e \cdot x ||_{L^{\infty}(B_1)} \leq \frac{\delta}{3}  .
\end{equation}
By the uniform convergence of \( (u_j) \) and (ii) we have
\begin{equation}
 \label{teoria3.2}
||u_r(x) - u_{j,r}(x)||_{L^{\infty}(B_1)} \leq \frac{\delta}{3} 
\end{equation}
and
\begin{equation}
 \label{teoria3.3}
||u_{j,r}(x) - v_{j,r}(x) ||_{L^{\infty}(B_1)} \leq \frac{\delta}{3}  
\end{equation} 
when \( j \) is large. Combining (\ref{teoria3.1}), (\ref{teoria3.2}) and (\ref{teoria3.3}) we obtain 
\[
 || u_{j,r} (x) - e \cdot x ||_{L^{\infty}(B_1)} \leq \delta
\]
and
\[
 || v_{j,r} (x) - e \cdot x ||_{L^{\infty}(B_1)} \leq \delta = \delta (\epsilon_0 /3, H_{\alpha})  .
\]
Therefore by Theorem \ref{wangyu1}
\[
 | D u_{j,r }(0) - e | = | D u_j(0) - e | \leq \epsilon_0 /3
\]
and
\[
 | D v_j(0) - e | \leq \epsilon_0 /3 \, .
\]
By triangle inequality \(| D u_j(0) - D v_j(0) | <  \epsilon_0 \), which contradicts \( (iii)\). \( \qquad \Box \)

\vspace*{4mm}

\textbf{Proof of Corollary 1.3.} We may assume that \(x_0 = 0 \). Since the other implication is trivial we only prove the nontrivial part of the corollary. For that define
\[
 k_0 = \inf \, \{ k > 0 \mid C_{k}^{H}(x) \geq 1 \, \text{for all} \,  x \in \partial \Omega \, \}  .
\]
Since \( u(y) = 1  \) and \( C_{k_0}^{H}(y) \geq 1 \) on \( \partial \Omega\) we have
\[
 \sup_{y,x \in \partial (\Omega \backslash \{0 \})} \{ u(y) - u(x) - C_{k_0}^{H}(y - x) \} \leq 0  .
\]
The general AMLE property yields
\begin{equation}
 \label{corollary0}
u(y) - u(x) - C_{k_0}^{H}(y - x)  \leq 0 
\end{equation}
for all \( y,x \in \overline{\Omega} \). In particular, 
\[
  u(y) \leq  C_{k_0}^{H}(y)
 \]
for all \(y \in \Omega \). We show next that   \( u(y) \geq  C_{k_0}^{H}(y) \) and we are done. 

Using the definition of \( k_0 \) we find \(\bar{x} \in \partial \Omega \) such that \( C_{k_0}^{H}(\bar{x}) = u(\hat{x}) = 1\). In this case it is easy to see that \( [0, \hat{x}] \subset \overline{\Omega} \). From this and from (\ref{corollary0}) it follows
\[
 u(t \bar{x}) = C_{k_0}^{H}(t \bar{x}) 
\]
for all \( t \in [0,1] \). Theorem \ref{teoria1} yields
\[
 \lim_{\lambda \to 0 +} \, \frac{u(\lambda x)}{ \lambda} =  C_{k_0}^{H}(x) 
\]
locally uniformly. Therefore by Theorem \ref{wangyu2} \( \lim_{x \to 0} | Du(x) - D C_{k_0}^{H}(x) | = 0\). Especially
\begin{equation}
\label{corollary 1}
 \lim_{x \to 0} H(Du(x))  = k_0 . 
\end{equation}

Choose \(y \in \Omega , \, y \neq 0\) such that \( Du(y) \neq 0 \) and a path   \(\xi  \)  such that \( \xi(0) = y \) and 
\[
\dot{\xi}(t) = - H_p(Du(\xi(t)))  . 
\]
Since \( u \in C^2 \) is a solution of Aronsson equation (\ref{aronsson2}) we have
\[ 
H(Du(\xi(t))) = \text{constant}  .
 \]
This implies that \( \xi \) can not stay inside \( \Omega \backslash \{0 \} \) forever, but there has to be \( \delta > 0 \)  such that  \( \xi(\delta) = 0 \). By (\ref{corollary 1}) \( \lim_{t \to\delta} \, H(Du(\xi(t))) = k_0 \) and therefore 
\begin{equation}
 \label{corollary2}
 H(Du(\xi(t))) = k_0  
\end{equation}
for all \( t \in [0, \delta]\). In particular, we have
\begin{equation}
\label{corollary3}
 H(Du(y)) = k_0 
\end{equation}
for all \( y \in  \Omega \backslash \{0 \}\) for which \( Du(y) \neq 0\).

From Proposition \ref{kartiolle2} \( (ii) \) we obtain that, if \( q \in \{ p \mid H(p) = k_0 \} \), then
\begin{equation}
 \label{corollary4}
C_{k_0}^{H}(H_p(q)) = q \cdot H_p(q)   .
\end{equation}
Hence
\[
\begin{split}
 u(y) = u(y) - u(0) &= \int\limits_{\delta}^0 Du(\xi(t)) \cdot \dot{\xi}(t) \, dt = \int\limits_{0}^{\delta} Du(\xi(t)) \cdot H_p(Du(\xi(t)))  \, dt \\
		&= \int\limits_{0}^{\delta} C_{k_0}^{H}[ H_p(Du(\xi(t))) ]  = \int\limits_{0}^{\delta} C_{k_0}^{H}(- \dot{\xi}(t))   \,  dt \geq C_{k_0}^{H}(- \int\limits_{0}^{\delta}  \dot{\xi}(t)) = C_{k_0}^{H}( y),
\end{split}	
\]
where the fourth equality follows from  (\ref{corollary2}) and (\ref{corollary4}) and the inequality is just Jensen's inequality. The use of Jensen's inequality is justified by the convexity and homogeneity of \(C_{k_0}^{H} \). Therefore 
\begin{equation}
\label{corollary5}
  u(y) \geq  C_{k_0}^{H}(y)  
\end{equation}
for all \( y \in \{ x \in \Omega \backslash \{ 0 \} \mid Du(x) \neq 0\} \eqqcolon V\). 

The argument (\ref{corollary3}) implies that \( V = \{ x \in \Omega \backslash \{ 0 \} \mid H(Du(x)) = k_0 \}\). Therefore \( V \) is both open and closed set. Hence \( V = \Omega \backslash \{ 0 \} \) and the inequality (\ref{corollary5}) holds for all \(y \in \Omega \).  \( \qquad \Box \)

\end{document}